\documentclass[12pt, reqno]{amsart}

\usepackage{amsmath,amsfonts,amsthm,amssymb,amscd}
\usepackage[utf8]{inputenc}
\usepackage[T2A]{fontenc}
\usepackage[all]{xy}
\usepackage{enumitem}
\usepackage[english]{babel}
\usepackage{graphicx}
\usepackage{xcolor}
\usepackage[margin=2.8cm]{geometry}

\usepackage{hyperref}
\hypersetup{
    hypertexnames=false,
    colorlinks,
    linkcolor={red!50!black},
    citecolor={blue!50!black},
    urlcolor={blue!80!black}
}

\newtheorem{theorem}[subsection]{Theorem}
\newtheorem{lemma}[subsection]{Lemma}

\newtheorem{definition}[subsection]{Definition}

\newtheorem{remark}[subsection]{Remark}

\makeatletter
\@addtoreset{subsection}{section}
\@addtoreset{equation}{section}
\@addtoreset{figure}{section}
\@addtoreset{table}{section}
\makeatother

\makeatletter
\newcommand\testshape{family=\f@family; series=\f@series; shape=\f@shape.}
\def\myemphInternal#1{\if n\f@shape%
\begingroup\itshape #1\endgroup\/%
\else\begingroup\sf\itshape #1\endgroup%
\fi}
\def\myemph{\futurelet\testchar\MaybeOptArgmyemph}
\def\MaybeOptArgmyemph{\ifx[\testchar \let\next\OptArgmyemph
                 \else \let\next\NoOptArgmyemph \fi \next}
\def\OptArgmyemph[#1]#2{\index{#1}\myemphInternal{#2}}
\def\NoOptArgmyemph#1{\myemphInternal{#1}}
\makeatother

\newcommand\id{\mathrm{id}}          

\newcommand\bR{\mathbb{R}}

\newcommand\Diff{\mathcal{D}}       
\newcommand\Homeo{\mathcal{H}}      
\newcommand\Orb{\mathcal{O}}        
\newcommand\Stab{\mathcal{S}}       

\newcommand\DiffId{\Diff_{\id}}     

\newcommand\Ci[2]{\mathcal{C}^{\infty}(#1,#2)}               
\newcommand\Morse[2]{\mathcal{M}(#1,#2)}                     

\newcommand\Stabilizer[1]{\Stab(#1)}

\newcommand\StabilizerIsotId[1]{\Stab'(#1)}

\newcommand\OrbitComp[2]{\Orb_{#1}(#2)}

\newcommand\Mman{M}

\newcommand\DiffM{\Diff(\Mman)}
\newcommand\DiffIdM{\DiffId(\Mman)}

\newcommand\func{f}

\newcommand\dif{h}

\newcommand\fa{f|_{A}}
\newcommand\fb{f|_{B}}
\newcommand\fG{G_{f}}
\newcommand\aG{G_{f|_{A}}}
\newcommand\bG{G_{f|_{B}}}
\newcommand\Sf{\mathcal{S}'({f})}
\newcommand\Sa{\mathcal{S}'({f|_{A}})}
\newcommand\Sb{\mathcal{S}'({f|_{B}})}
\newcommand\Rf{\rho(\mathcal{S}'({f}))}

\newcommand\ag{\gamma|_{\gA}}
\newcommand\bg{\gamma|_{\gB}}

\newcommand\ha{\dif|_{A}}

\newcommand\rh{\rho(\dif)}
\newcommand\rha{\rho(\dif|_{A})}
\newcommand\rhb{\rho(\dif|_{B})}
\newcommand\gG{\varGamma}
\newcommand\gA{\varGamma_{A}}
\newcommand\gB{\varGamma_{B}}

\newcommand\KRGraph[1]{\varGamma_{#1}}
\newcommand\KRGraphf{\KRGraph{\func}}
\newcommand\FSp[1]{\mathcal{F}(\Mman,\bR)}

\author{Anna Kravchenko}
\email{annakravchenko1606@gmail.com}
\address{Department of Geometry, Topology, and Dynamical Systems \\
Taras Shevchenko National University of Kyiv
Hlushkova Avenue, 4e, Kyiv, Ukraine, 03127}

\author{Sergiy Maksymenko}
\email{maks@imath.kiev.ua}
\address{Topology Laboratory of Algebra and Topology Department \\
 Institute of Mathematics of National Academy of Sciences of Ukraine \\
 Tereshchenkivs'ka str. 3, Kyiv, Ukraine, 01024}

\title[Automorphisms of Kronrod-Reeb graphs]{Automorphisms of Kronrod-Reeb graphs of Morse functions on $2$-sphere}


\begin{document}

\begin{abstract}
Let $M$ be a compact two-dimensional manifold, $f \in C^{\infty}(M,\mathbb{R})$ be a Morse function, and $\Gamma_f$ be its Kronrod-Reeb graph.
Denote by $\mathcal{O}(f)=\{f \circ h \mid h \in \mathcal{D}\}$ the orbit of $f$ with respect to the natural right action of the group of diffeomorphisms $\mathcal{D}$ on $C^{\infty}(M,\mathbb{R})$, and by $\mathcal{S}(f)=\{h\in\mathcal{D} \mid f \circ h = f\}$ the corresponding stabilizer of this function.
It is easy to show that each $h\in\mathcal{S}(f)$ induces a homeomorphism of $\Gamma_f$.
Let also $\mathcal{D}_{\mathrm{id}}(M)$ be the identity path component of $\mathcal{D}(M)$, $\mathcal{S}'(f)= \mathcal{S}(f) \cap \mathcal{D}_{\mathrm{id}}(M)$ be group of diffeomorphisms of $M$ preserving $f$ and isotopic to identity map, and $G_f$ be the group of homeomorphisms of the graph $\Gamma_f$ induced by diffeomorphisms belonging to $\mathcal{S}'(f)$.
This group is one of the key ingredients for calculating the homotopy type of the orbit $\mathcal{O}(f)$.

Recently the authors described the structure of groups $G_f$ for Morse functions on all orientable surfaces distinct from $2$-torus $T^2$ and $2$-sphere $S^2$.
The present paper is devoted to the case $\Mman=S^{2}$.
In this situation $\Gamma_f$ is always a tree, and therefore all elements of the group $G_f$ have a common fixed subtree $\mathrm{Fix}(G_f)$, which may even consist of a unique vertex.
Our main result calculates the groups $G_f$ for all Morse functions $\func\colon S^{2}\to\mathbb{R}$ whose fixed subtree $\mathrm{Fix}(G_f)$ consists of more than one point.
\end{abstract}


\maketitle

\section{Introduction}
Let $\Mman$ be a compact two-dimensional manifold and $\DiffM$ the group of diffeomorphisms of $\Mman$.
Then there exists a natural right action \[\phi\colon \Ci{\Mman}{\bR}\times\DiffM\to \Ci{\Mman}{\bR}\] of this group on the space of smooth functions on $\Mman$ defined by the formula $\phi(\func,\dif) = \func \circ \dif$.
For $\func \in \Ci{\Mman}{\bR}$ denote by
\begin{align*}
\Stabilizer{\func} &=\{\dif\in\DiffM \mid \func \circ \dif = \func\}
\end{align*}
its \myemph {stabilizer} with respect to the specified action.

\begin{definition}\label{def:MorseFunc}
Let $\FSp{\Mman}$ be the subset of $\Ci{\Mman}{\bR}$ consisting of maps $\func\colon \Mman\to\bR$ such 
\begin{enumerate}[leftmargin=*, label={\rm(\arabic*)}]
\item $\func$ takes constant values on the connected components of the boundary $\partial\Mman$ and has no critical points on $\partial\Mman$;
\item for each critical point $z$ of $\func$ there are local coordinates $(x,y)$ in which $z=(0,0)$ and $\func(x,y)=f(z) + g_z(x,y)$, where $g_z\colon \bR^2\to\bR$ is a homogeneous polynomial without multiple factors.
\end{enumerate}
Notice that every critical point of $\func\in\FSp{\Mman}$ is isolated.

A function $\func\in\FSp{\Mman}$ is called \myemph{Morse}, if $\deg g_z =2$ for each critical point $z$ of $\func$.
In that case, due to Morse Lemma, one can assume that $g_z(x,y) = \pm x^2 \pm y^2$.
\end{definition}    

We will denote by $\Morse{\Mman}{\bR}$ the space of all Morse maps $\Mman\to\bR$.

Homotopy types of stabilizers and orbits of Morse functions and functions from $\mathcal{F}(\Mman,\bR)$ were studied in \cite{Maksymenko:AGAG:2006}, \cite{Maksymenko:ProcIM:ENG:2010}, \cite{Maksymenko:UMZ:ENG:2012}, \cite{Kudryavtseva:ConComp:VMU:2012}, \cite{Kudryavtseva:MathNotes:2012}, \cite{Kudryavtseva:MatSb:2013}, \cite{KudryavtsevaPermyakov:MatSb:2010}, \cite{MaksymenkoFeshchenko:MS:2015}, \cite{MaksymenkoFeshchenko:MFAT:2015}.

Let $\func\in\Ci{\Mman}{\bR}$, $\KRGraphf$ be a partition of the surface $\Mman$ into the connected components of level sets of this function, and  $p\colon \Mman \to \KRGraphf$ be the canonical factor-mapping, associating to each $x \in \Mman$ the connected component of the level set $\func^{-1}(\func(x))$ containing that point.

Endow $\KRGraphf$ with the factor topology with respect to the mapping $p$: so a subset $A\subset \KRGraphf$ will be regarded as open if and only if its inverse image $p^{-1}(A)$ is open in $\Mman$.
Then $\func$ induces the function $\hat{\func}\colon \KRGraphf \to \bR$, such that $\func=\hat{\func}\circ p$.

It is well known, that if $\func\in\FSp{\Mman}$, then $\KRGraphf$ has a structure of a one-dimensional CW-complex called the \myemph{Kronrod-Reeb graph}, or simply the \myemph{graph} of $\func$.
The vertices of this graph correspond to critical connected components of level sets of $\func$ and connected components of the boundary of the surface.
By the \myemph{edge} of $\KRGraphf$ we will mean an \myemph{open} edge, that is, a one-dimensional cell.

Denote by $\Homeo(\KRGraphf)$ the group of homeomorphisms of $\KRGraphf$.
Notice that each element of the stabilizer $\dif\in\Stabilizer{\func}$ leaves invariant each level set of $\func$, and therefore induces a homeomorphism $\rho(\dif)$ of the graph of $\func$, so that the following diagram is commutative:
\begin{equation}\label{equ:2x2_M_Graph}
\xymatrix{
	\Mman \ar[rr]^-{p} \ar[d]_-{\dif} &&
	\KRGraphf \ar[rr]^-{\hat{\func}} \ar[d]^-{\rho(\dif)} &&
	\bR \ar@{=}[d]  \\
	\Mman \ar[rr]^-{p} &&
	\KRGraphf \ar[rr]^-{\hat{\func}} &&
	\bR
}
\end{equation}

Moreover, the correspondence $h\mapsto g(h)$ is a homomorphism of groups \[\rho\colon \Stabilizer{\func} \to \Homeo(\KRGraphf).\]

Let also $\DiffIdM$ be the path component of the identity map $\id_{\Mman}$ in $\DiffM$. 
Put 
\begin{align*}
\StabilizerIsotId{\func} &= \Stabilizer{\func} \cap \DiffIdM   &\
\fG &=\rho(\StabilizerIsotId{\func}).
\end{align*}

Thus, $\fG$ is the group of automorphisms of the Kronrod-Reeb graph of $\func$ induced by diffeomorphisms of the surface preserving the function and isotopic identity.

\begin{remark}\label{rm:f}\rm
Since $\hat{\func}$ is monotone on edges of  $\KRGraphf$, it is easy to show that $\fG$ is a finite group.
Moreover, if $g(E)=E$, for some $g\in G$ and an edge $E$ of the graph $\KRGraphf$, then $g(x)=x$ for all $x \in E$.
\end{remark}

Since $\fG$ is finite and $\rho$ is continuous, it follows that $\rho$ reduces to an epimorphism
\[
 \rho_0\colon \pi_0 \StabilizerIsotId{\func} \to \fG,
\]
of the group $\pi_{0}\StabilizerIsotId{\func}$ path components of $\StabilizerIsotId{\func}$ being an analogue of the mapping class group for $\func$-preserving diffeomorphisms.

Algebraic structure of the group $\pi_{0}\StabilizerIsotId{\func}$ of connected components of $\StabilizerIsotId{\func}$ for all $\func\in\FSp{\Mman}$ on orientable surfaces $\Mman$ distinct from $2$-torus and $2$-sphere is described in~\cite{Maksymenko:KRGraphs:2013}, and the structure of its factor group $\fG$ is investigated in~\cite{MaksymenkoKravchenko:GMF:2018}.
These groups play an important role in computing the homotopy type of the path component $\OrbitComp{\func}{\func}$ of the orbit of $\func$, see also~\cite{Maksymenko:AGAG:2006}, \cite{Maksymenko:ProcIM:ENG:2010}, \cite{Kudryavtseva:ConComp:VMU:2012}, \cite{Kudryavtseva:MathNotes:2012}, \cite{Kudryavtseva:MatSb:2013}.

The purpose of this note is to describe the groups $\fG$ for a certain class of smooth functions on $2$-sphere $S^2$.

The main result Theorem~\ref{th:iso} reduces computation of $\fG$ to computations of similar groups for restrictions of $\func$ to some disks in $S^2$.
As noted above the latter calculations were described in \cite{MaksymenkoKravchenko:GMF:2018}.

First we recall a variant of the well known fact about automorphisms of finite trees from graphs theory.

\begin{lemma}\label{lm:cw}
Let $\Gamma$ be a finite contractible one-dimensional CW-complex (<<a topological tree>>), $G$ be a finite group of its cellular homeomorphisms, and $\mathrm{Fix}(G)$ be the set of common fixed  points of all elements of the group $G$.
Then $\mathrm{Fix}(G)$ is either a contractible subcomplex or consists of a single point belonging to some edge $E$ an open 1-cell), and in the latter case there exists $g \in G$ such that $g(E)=E$ and $g$ changes the orientation of $E$.
\end{lemma}

Suppose $\func\colon \mathcal{S}^{2}\to{\bR}$ belongs to $\FSp{S^2}$.
Then it is easy to show that $\KRGraphf$ is a tree, \emph{i.e.}, a finite contractible one-dimensional CW-complex, and by Remark~\ref{rm:f} $\fG$ is a finite group of cellular homeomorphisms of $\KRGraphf$.
Therefore, for $\fG$, the conditions of Lemma~\ref{lm:cw} are satisfied.
Note that according to Remark~\ref{rm:f} the second case of Lemma~\ref{lm:cw} is impossible, and hence $\fG$ has a fixed subtree.

In this paper we consider the case when the fixed subtree of the group $\fG$ contains more than one vertex, \emph{i.e.} has at least one edge.

Let us also mention that $\DiffId(S^{2})$ coincides with the group $\Diff^{+}(S^{2})$ of diffeomorphisms of the sphere preserving orientation, \cite{Smale:ProcAMS:1959}.
Therefore $\StabilizerIsotId{\func}$ consists of diffeomorphisms of the sphere preserving the function $f$ and the orientation of $\mathcal{S}^{2}$.
\begin{theorem}\label{th:iso}
Let $\func\in\FSp{S^2}$.
Suppose that all elements of the group $\fG$ have a common fixed edge $E$.
Let $x \in E$ be an arbitrary point and $A$ and $B$ be the closures of the connected components of $S^{2}\setminus p^{-1}(x)$.
Then
\begin{enumerate}[label={\rm{(\arabic*)}}]
\item\label{en1} $A$ and $B$ are 2-disks being invariant with respect to $\StabilizerIsotId{\func}$;
\item\label{en2}
the restrictions $\fa \in\FSp{A}$ and $\fb\in\FSp{B}$;
\item\label{en3}
the map $\phi\colon \fG \to \aG \times \bG$ defined by the formula 
\[\phi(\gamma)=(\ag, \bg)\] is an isomorphism of groups.
\end{enumerate}

\begin{proof}
\ref{en1}
By assumption $x$ belongs to the open edge $E$.
Therefore $p^{-1}(x)$ is a regular connected component of some level set of the function $\func$, that is, a simple closed curve.
Then, by Jordan Theorem, $p^{-1}(x)$ divides the sphere into two connected components whose closures are homeomorphic to two-dimensional disks.
Consequently, $A$ and $B$ are two-dimensional disks.

Let as show that $A$ and $B$ are invariant with respect to $\StabilizerIsotId{\func}$, \emph{i.e.}, $h(A)=A$ and $h(B)=B$ for each $h\in \StabilizerIsotId{\func}$.
Denote
\begin{align*}
\gA&=p(A) &\
\gB&=p(B).
\end{align*}
Then 
\begin{align*}\label{equ:gG}
\gA\cup\gB&=\gG &\
\gA\cap\gB&=\{x\}.
\end{align*}
By definition, $\rh(x)=x$, whence $\rh$ either preserves both $\gA$ and $\gB$ or interchange them.
We claim that
\begin{align*}
\rh(\gA)&=\gA &\
\rh(\gB)&=\gB.
\end{align*}
Indeed suppose $\rh(\gA)=\gB$.
Since $\rh$ is fixed on $E$, it follows that
\[
\rh(\gA\cap E)=\gA\cap E,
\]
whence
\[
\rho(h)(\gA\cap E)=\rh(\gA)\cap\rho(E)=\gB\cap E\neq\gA\cap E,
\]
which contradicts to our assumption.
Thus $\gA$ and $\gB$ are invariant with respect to the group $\fG$.

\medskip 

Now we can show that $A$ and $B$ are also invariant with respect to $h$.
By virtue of the commutativity of the diagram~\eqref{equ:2x2_M_Graph} $\rh(p(y))=p(h(y))$ for all $y\in \gG$.
In particular:
\[
p(h(A))=\rh(p(A))=\rh(\gA)=\gA.
\]

Therefore, $h(A)=p^{-1}(\gA)=A$.
The proof for $B$ is similar.
Thus, $A$ and $B$ are invariant with respect to $\StabilizerIsotId{\func}$.

\medskip 

\ref{en2} 
Notice that the function $f$ takes a constant value on the simple closed curve $p^{-1}(x)$ being a common boundary of disks $A$ and $B$, and does not contain critical points of $\func$.
Therefore, the restrictions $\fa, \fb$ satisfy the conditions 1) and 2) the Definition~\ref{def:MorseFunc}, and so they belong to 
$\FSp{A}$ and $\FSp{B}$ respectively.

\medskip 

\ref{en3}
We should prove that the map $\phi\colon \fG \to \aG \times \bG$ defined by formula  $\phi(\gamma)=(\ag,\bg)$ is an isomorphism. 

First we will show that $\phi$ is correctly defined.
Let $\gamma\in\fG=\Rf$, that is, $\gamma=\rh$, where $\dif$ is a diffeomorphism of the sphere preserving the function $f$ and isotopic to the identity.

We claim that $\ha\in\Sa=\Stabilizer{\fa}\cap \DiffId(A)$.
Indeed, for each point $x\in A$ we have that:
\[
f(x)=\fa(x)=\fa(\ha(x))=\fa(\dif(x))=f(\dif(x)),
\]
which means that $\ha\in\Stabilizer{\fa}$.

Moreover, since $\dif$ preserves the orientation of the sphere, it follows that $\ha$ preserves the orientation of the disk $A$, and therefore by~\cite{Smale:ProcAMS:1959}, $\ha\in\DiffId(A)$.
Thus $\ag\in\aG$.
Similarly $\bg\in\bG$, and so $\phi$ is well defined.

Let us now verify that $\phi$ is an \myemph{isomorphism of groups}, that is, a bijective homomorphism.

Let $\delta, \omega \in\fG$.
Then
\begin{align*}
\phi(\delta\circ\omega)
&= \bigl( \delta\circ\omega|_{\gA}, \ \delta\circ\omega|_{\gB} \bigr) = \\
&=\bigl(\delta|_{\gA},\ \delta|_{\gB}\bigr)\circ\bigl(\omega|_{\gA}, \ \omega|_{\gB}\bigr)=\\
&=\bigl(\delta|_{\gA}\circ\omega|_{\gA},\ \delta|_{\gB}\circ\omega|_{\gB}\bigr)=\\
&=\bigl(\delta\circ\omega|_{\gA}, \ \delta\circ \omega|_{\gB}\bigr),
\end{align*}
sp $\phi$ is a homomorphism.

Let us show that $\ker\phi=\{\id_{\gG}\}$.
Indeed, suppose $\gamma\in \ker\phi$, that is $\ag=\id_{\gA}$ and $\bg=\id_{\gB}$.
Then $\gamma$ is fixed on $\gA\cup\gB=\gG$, and hence it is the identity map.

Surjectivity of $\phi \colon \fG \to \aG \times \bG$ is implied by the following simple lemma whose proof we leave to the reader.

\begin{lemma}\label{lm:ab}
Suppose $f\colon D^{2}\to \bR$ belgns to the space $\FSp{D^{2}}$.
Then for arbitrary $\alpha\in\fG$, there exists $a\in \Sf$ fixed near the boundary $\partial D^{2}$ and such that $\alpha=\rho(a)$.\qed
\end{lemma}

Let $(\alpha,\beta)\in\aG \times \bG$, then by Lemma~\ref{lm:ab} there exist $a\in \Sa$ and $b\in \Sb$ fixed near $\partial A=\partial B=p^{-1}(x)$ and such that $\alpha=\rho_{A}(a)$ and $\beta=\rho_{B}(b)$.
Define $\dif$ by the following formula:
\begin{equation*}
   \dif=
   \begin{cases}
         a(x),&  x\in A, \\
         b(x),&  x\in B.
   \end{cases}
\end{equation*}
Then, $\dif$ is a diffeomorphism of the sphere, preserving the function and orientation, whence $\dif\in\StabilizerIsotId{\func}$.

Moreover if we put $\gamma=\rh\in \fG$, then $\ag=\rha=\alpha$ and $\bg=\rhb=\beta$.
In other words, $\phi(\gamma)=(\ag,\bg)=(\alpha,\beta)$, \emph{i.e.}, $\phi$ is surjective and therefore an isomorphism.
\end{proof}
\end{theorem}



\def\cprime{$'$}

\end{document}